\documentclass[11pt, a4paper]{article}
\usepackage{amsmath}
\usepackage{changepage}
\usepackage{enumitem}
\usepackage{hyperref}

\title{Undergraduate algebra in nineteenth-century Oxford}
\author{\textsc{Joseph Gage} \\
\href{mailto:joseph.gage@lmh.ox.ac.uk}{joseph.gage@lmh.ox.ac.uk} \\ 
\textit{University of Oxford}}
\date{}	

\begin{document}

\maketitle

\begin{abstract}
The nineteenth century was an important period for both Oxford mathematics and algebra in general. While there is extensive documentation of mathematical research in Oxford at this time, the same cannot be said of the teaching. The content of the course presents a different picture: it shows what those who set it felt was most valuable for a young mathematician to learn, perhaps indicating what direction they expected mathematics to take in the future. To find out what undergraduates were taught, I have looked through examination papers of the years between 1828 and 1912 with a focus on algebra, as well as supporting material. In this paper I will present my findings. I will give a picture of what an Oxford undergraduate's course in algebra looked like by the end of the nineteenth century and discuss my own conclusions as to why it took such a form.
\end{abstract}

\noindent
Algebra's development in the nineteenth century is well documented. Due to a number of factors including the Leibniz--Newton priority dispute, relations between the respective mathematical communities of the UK and the rest of Europe had deteriorated, leaving English universities isolated from much of the progress made on the continent (Katz 1993, 482). The mainland mathematicians of this time can perhaps claim primary responsibility for the growth of algebra into the more abstract subject it is today, and much research went into this rapidly evolving field. In England, efforts in this direction by the likes of Cayley and De Morgan were met with scepticism by their contemporaries. Mathematical physicist William Thomson lamented in a letter of 1864 that they would `devote what skill they had \ldots\ to pieces of algebra which possibly interest four people in the world' (Thompson 1910, 433).

The view of algebra in Oxford during this period is muddled slightly by the low status of mathematics there at the time. Oxford could not match the prestige held by the Cambridge tripos, and there was a feeling among some in the university that mathematical study should not be a priority for undergraduates. In 1817, Charles Atmore Ogilvie, fellow of Balliol College, wrote in a prize essay that `by such an union of classical with mathematical learning as makes the latter in an important matter subsidiary to the former, the student will be best qualified for active life' (Ogilvie 1817, 191). For a brief period between 1800 and 1807 mathematical study had been compulsory for anyone wishing to graduate, but this did not last and mathematics remained an optional course to be taken alongside the Literae Humaniores BA until 1864 (Hannabuss 1997, 444).

Nevertheless, much can be learned from the mathematical exam papers of this period. With drastic changes taking place in what was considered \textit{algebra} it is enlightening to look into how these changes made their way into the taught course at Oxford. Hannabuss (1997, 450) credits Henry Smith, Savilian Professor of Geometry\footnote{One of the three mathematical chairs in place at Oxford at this point, alongside the Savilian Professorship of Astronomy and the Sedleian Professorship of Natural Philosophy. These remained the only chairs until the creation of the Waynflete Professorship of Pure Mathematics in 1892 (Busbridge 1974).} 1861--83, with bringing Oxford mathematics `international renown' during his tenure. Perhaps then at last beginning to feel the influence of continental ideas, a distinct flavour of algebra becomes apparent in the undergraduate course material, and can be seen in the papers set to students each year.

\section*{The mathematical examination}

Prior to the nineteenth century there was not much in the way of a rigorous examination system in place at Oxford. Candidates were required to partake in three disputations and one oral exam, but these were cursory and failure was almost unheard-of (Sutherland 1986, 475). When initial reform came in 1800 the oral exam was retained, but a pass could no longer be assured: the failure rate increased from less than five per cent to around twenty per cent after 1810 and there are reports of both students and examiners collapsing from the stress of exams (Curthoys 1997, 346).

The written exam emerged out of the questions set for candidates while they waited to be examined orally (Roberts 1814, 73). With the number of students seeking to graduate increasing from 188 in 1810 to 404 in 1827, the amount of written work set increased until, in 1828, the first printed papers in something like the modern style were set for mathematics candidates. However, these did not immediately replace the viva, as this format was preferred for the tests on religion and scripture that every candidate was required to undertake until 1864\footnote{It was not until 1864 that single-subject degrees were permitted in any subject besides classics. Up to that point, candidates hoping to study mathematics (or, from 1850, natural science or law and modern history) would be required to do so in addition to their classical studies. They would then be examined in both subjects---attaining first class honours in both was the original meaning of the `double first' (Curthoys 1997, 352).} (Curthoys 1997, 348).

Moderations, or the `first public examination', were introduced in 1850. Accordingly, the final examination then became the `second public examination'. Previously undergraduates would only undergo Responsions before their final exam---these were simple tests in place of an entrance exam and were taken early in the degree. Moderations took place at the end of a student's second year and were considerably more testing (Curthoys 1997, 354). In terms of content, the first and second public examinations in mathematics were similar. The difference came in the increased difficulty for the final examination, suggesting a more linear course than the modular structure in place today.

The examination initially consisted of ten papers: Algebra, Geometry, etc. (three papers), Differential and Integral Calculus (two papers), Astronomy, Hydrostatics, (Newton's) Principia\footnote{The Oxford  course maintained its obsession with Newton throughout the nineteenth century. As well as having an entire paper dedicated to his \textit{Philosophi\ae\ Naturalis Principia Mathematica}, his mathematics pervaded other papers. Early papers would explicitly ask candidates to describe his techniques. The applied courses in optics, mechanics and astronomy were unsurprisingly informed by his ideas, but his name makes regular appearances in the pure papers too. Candidates were regularly asked to reproduce his proof of the generalized binomial theorem, and in the sections on the theory of equations his method for approximation of roots was a fixture for the latter half of the century.}, Optics, and Mechanics. Each paper would contain between six and nine questions---it is not clear whether candidates were expected to answer all questions or some subset thereof. When the mathematical scholarship was instated by Baden Powell in 1831 (Hannabuss 1997, 446), it was awarded based on five papers of around 12 questions each. The format and content of the exams were in a state of constant, if slow, change but the fundamentals were consistent: a roughly even split between pure and applied mathematics, with a general problems paper that drew on material from multiple courses later introduced. As well as algebra, the pure mathematics papers dealt largely with geometry, trigonometry, and calculus; basic questions on real analysis, combinatorics, and probability theory were also common.

\section*{Solving equations}

It is unlikely that the first algebra questions would inspire the same fear in students today as they did in the early nineteenth century: during the 1830s the material examined rarely extended further than what would be expected of sixth form students today. In fact, for the first thirty or so years of examinations there was little change in the topics tested. There was a steady increase in difficulty of the questions, however, illustrated by the following similar questions, dated fifteen years apart:

\vspace{3mm}
\begin{adjustwidth}{1cm}{1cm}
\textit{4. Solve the following equations:}
\begin{enumerate}[label=(\arabic*), topsep=2mm, itemsep=1mm]
\item $\sqrt{x} - 4 = \dfrac{259 - 10x}{\sqrt{x} + 4}$.
\item $3x + \dfrac{3}{x} = 10$.
\item $x^{3} + 3x - 14 = 0$.
\end{enumerate}
\end{adjustwidth}
\begin{adjustwidth}{0.5cm}{0.5cm}
\begin{flushright}
(Mathematical Examination: Easter Term, 1839)
\end{flushright}
\end{adjustwidth}
\vspace{3mm}
\begin{adjustwidth}{1cm}{1cm}
\textit{11. Solve the following equations:}
\begin{enumerate}[label=(\arabic*), topsep=2mm, itemsep=1mm]
\item $x^{3} + 2x + 12 = 0$.
\item $x^{3} - 9x + 9 = 0$.
\item $x^{10} - 5x^{8} + 9x^{6} + 9x^{4} + 5x^{2} + 1 = 0$.
\item $x^{4} + x^{3} - 8x^{2} - 16x - 8 = 0$, where $1 - \sqrt{5}$ is a root.
\item $x^{4} - 6x^{3} + 14x^{2} - 15x + 6 = 0$, where there are two pairs of roots having their sums respectively equal.
\end{enumerate}
\end{adjustwidth}
\begin{adjustwidth}{0.5cm}{0.5cm}
\begin{flushright}
(Mathematical Examination: Easter Term, 1854)
\end{flushright}
\end{adjustwidth}
\vspace{3mm}

This style of question was typical of the algebra papers for much of the nineteenth century. It seems that students were taught various methods for solving specific classes of polynomial equations (usually up to degree five at most though occasionally higher), which would be called on in the examination. A popular choice was Cardano's method for the solution of a cubic of the form $x^{3} + px + q = 0$. Occasionally information would be given about the roots of the equation, as in the second example above, perhaps indicating that these methods were more crude for higher degree polynomials.

In general, the theory of equations, especially polynomials, was a very popular topic. Proofs of the rational root theorem, the complex conjugate root theorem, and Descartes' rule of signs were commonly asked for. As well as simply solving equations, candidates would frequently be asked about the relationships between a polynomial's roots and its coefficients. Knowledge of Vieta's formulae and Newton's identities was presumed and often tested\footnote{Vieta's formulae directly relate a polynomial's coefficients to the elementary symmetric polynomials in its roots. Newton's identites (or the Newton--Girard formulae) are a series of expressions for the sums of powers of $n$ variables in terms of the same elementary symmetric polynomials in those variables. They can be applied together in order to find a sum of powers of a polynomial's roots given its coefficients, which was a common exam question.}. Historically, these were an early appearance of symmetric polynomials, which were an important concept in nineteenth-century algebra---particularly in Galois theory. Despite the fact that texts including Galois theory had been in existence since 1866 (Katz 1993, 603), the first explicit treatment of symmetric polynomials in the Oxford examinations did not come until much later:

\vspace{3mm}
\begin{adjustwidth}{1cm}{1cm}
\textit{2. Prove that any rational integral symmetric function of $n$ variables $a_{1}, a_{2}, \ldots , a_{n}$ can be expressed in terms of the fundamental functions
\begin{align*}
\sum a_{1}, \sum a_{1}a_{2}, \sum a_{1}a_{2}a_{3}, \ldots , a_{1}a_{2} \ldots a_{n}.
\end{align*}}

\textit{Calling these latter $p_{1}, p_{2}, \ldots , p_{n}$ express $\sum a_{1}^{2}a_{2}^{2}a_{3}^{2}$ in terms of them.}
\end{adjustwidth}
\begin{adjustwidth}{0.5cm}{0.5cm}
\begin{flushright}
(Second Public Examination: Trinity Term, 1901)
\end{flushright}
\end{adjustwidth}
\vspace{3mm}

This is the familiar fundamental theorem of symmetric polynomials. While the result dates back to the eighteenth century (O Neumann 2007, 108), in the greater context that had been attached to symmetric polynomials by this point this may be the first sign of something resembling Galois theory in the Oxford course. There is, however, no further indication of this in the papers of the years immediately following, nor any mention of permutations---the key idea behind symmetric polynomials' relevance.

Another popular style of question was transforming polynomials. Most often candidates would be given a general polynomial and its roots $\alpha , \beta , \gamma , \ldots\ $ and be required to construct a polynomial with roots given as specific functions of the original roots, for example $\alpha + \beta , \beta + \gamma , \ldots\ $. An alternative would be to transform a polynomial by a change of variables into another that lacked a certain term, for example the square term in a cubic---transforming a general cubic into one that can be solved by Cardano's method. This was the first foray into the subject of linear transformations on polynomials (that is, a map of the form $x \mapsto ax + b$), which is something we will revisit later.

\section*{Letters and numbers}

Algebra has long concerned symbols and the rules for their manipulation. Indeed, to a typical modern secondary school student in the UK, algebra \textit{is} letters and symbols. Such a student would recognize the style of many questions from the nineteenth-century Oxford papers. From the introduction of the written exam until the final quarter of the century, especially in the first public examination papers of the 1850s and 60s, there was an overwhelming presence of questions that tested a candidate's ability to manipulate numbers and symbols. Many were classic number puzzles, sometimes entirely symbolic:

\vspace{3mm}
\begin{adjustwidth}{1cm}{1cm}
\textit{2. If $a, b, c, \ldots$ are the numbers of days in which $A, B, C, \ldots$ could respectively finish a piece of work alone; find how many days it will take them when working all together.}
\end{adjustwidth}
\begin{adjustwidth}{0.5cm}{0.5cm}
\begin{flushright}
(Mathematical Examination: Easter Term, 1842)
\end{flushright}
\end{adjustwidth}
\vspace{3mm}

Another familiar problem that was regularly set concerns two coaches travelling from London to York and York to London respectively, crossing over at a given time and place---the aim being to find the times of arrival. Questions were occasionally set in a `real-world' context; financial terms were especially popular in the middle years of the century:

\vspace{3mm}
\begin{adjustwidth}{1cm}{1cm}
\textit{1. $A$ lends $B$ a sum of money, which $B$ undertakes to repay with interest by $n$ annual instalments, the interest for each year being calculated on the portion of the debt left unpaid at the beginning of the year. Find the amount of each yearly payment.}
\end{adjustwidth}
\begin{adjustwidth}{0.5cm}{0.5cm}
\begin{flushright}
(First Public Examination: Michaelmas Term, 1865)
\end{flushright}
\end{adjustwidth}
\vspace{3mm}

This would suggest an inclination of the mathematics course to prepare students for life outside of academia. However, if this was the attitude at the time then it had changed by the end of the century as such questions were phased out. Symbolic manipulation was by this point a tool used in questions on other topics, a prime example being the theory of determinants.

The first questions on determinants appeared in the 1870s. Both proofs of properties of determinants (multilinearity, multiplicativity, and the alternating property) and explicit computation were required. The study of determinants emerged separately to that of matrices (Katz 1993, 621), and the concept of a matrix as something beyond a determinant was not something that ever appeared in the Oxford exams in the nineteenth century. It does seem, however, that the teaching of determinants did show shades of the wider theory of matrices. The multiplicative rule for determinants is the same as that of matrices, and the exam papers presented determinants in the form of an array. There was also the following question, which hinted at eigenvalues and the first part of the spectral theorem for real symmetric matrices:

\vspace{3mm}
\begin{adjustwidth}{1cm}{1cm}
\textit{4. Establish the rule for the multiplication of two determinants, and hence by forming the product $f(\lambda)f(-\lambda)$, shew that all the roots of $f(\lambda) = 0$ are real, where
\begin{align*}
f(\lambda) \equiv \begin{vmatrix}
a - \lambda, & h, & g \\ 
h, & b - \lambda, & f \\ 
g, & f, & c - \lambda
\end{vmatrix}.
\end{align*}
}
\end{adjustwidth}
\begin{adjustwidth}{0.5cm}{0.5cm}
\begin{flushright}
(Second Public Examination: Trinity Term, 1887)
\end{flushright}
\end{adjustwidth}
\vspace{3mm}

Rather than signifying the advent of linear algebra in the mathematics course, though, this appears to have been an anomaly. There were no similar questions in later years; determinants largely remained something studied for their own sake.

\section*{More on linear transformations}

One area where determinants did find an application in the exam papers was with reference to linear transformations. We have already seen that the idea of linear transformations on polynomials had been present in the Oxford course since the early years of exams, but the subject started to become more sophisticated with the introduction of invariants. The notion of an invariant is central to algebra. When transforming objects in a particular way, it is natural to look at what doesn't change---what are the fundamental properties that define the object? The outcome of this inquiry about linear transformations on polynomials was what is now known as \textit{invariant theory}.

It is perhaps unsurprising that invariant theory became one of the most popular aspects of algebra examined at Oxford by the end of the nineteenth century. First introduced by Cayley\footnote{Hermann Weyl (1939, 489) wrote that invariant theory `came into existence about the middle of the nineteenth century somewhat like Minerva: a grown-up virgin, mailed in the shining armor of algebra, she sprang forth from Cayley's Jovian head.'} in his 1845 paper `On the theory of linear transformations', one of its leading figures was James Joseph Sylvester, Smith's successor to the Geometry chair at Oxford\footnote{For a brief, but more detailed treatment of the emergence of invariant theory, see Parshall (2011): section entitled `The evolution of the theory of invariants'.}.

Rather than polynomials, teaching on invariant theory at Oxford generally focused on binary forms. Exam questions began appearing from the early 1880s onwards, and often began by asking candidates to define various terms from the theory---these were still relatively new concepts after all, and a departure from much of the familiar course material. The questions made links with partial differential equations, and determinants were clearly involved (in particular the Hessian):

\vspace{3mm}
\begin{adjustwidth}{1cm}{1cm}
\textit{9. Define the terms Invariant and Covariant of a Quantic, and prove that any invariant of a covariant is an invariant of the original quantic.}

\textit{If a binary quantic has any root repeated shew that the same root will occur repeated in its Hessian, and hence obtain the discriminant of the binary cubic $(a, b, c, d)(x, y)^{3}$ in the form $(ad-bc)^{2} - 4(ac - b^{2})(bd - c^{2})$.}
\end{adjustwidth}
\begin{adjustwidth}{0.5cm}{0.5cm}
\begin{flushright}
(Second Public Examination: Trinity Term, 1882)\footnote{The term `quantic' was more commonly used for algebraic forms at this time. The notation $(a_{0}, a_{1}, \ldots , a_{n})(x, y)^{n}$ gives the coefficients of the terms in the binary form. In this example, $(a, b, c, d)(x, y)^{3} = ax^{3} + bx^{2}y + cxy^{2} + dy^{3}$.}
\end{flushright}
\end{adjustwidth}
\vspace{3mm}

By the end of the century, invariant theory was a major part of the mathematics course at Oxford. Exam questions indicate that candidates were expected to be more familiar with the ideas behind them. However, the content was still very much non-abstract. Invariant theory today more generally concerns the action of groups on algebraic objects, but the Oxford course remained rooted in binary forms. The focus was by this point on various canonical forms for these forms and how invariants dictated them:

\vspace{3mm}
\begin{adjustwidth}{1cm}{1cm}
\textit{9. Prove that a binary cubic can be reduced by a linear transformation to the sum of two cubes.}

\textit{Prove also that a binary quartic can be reduced similarly to the form $X^{4} + Y^{4}$ only when its invariant $J$ (of order $3$) vanishes; and show how, when this is the case, to determine $X$ and $Y$.}
\end{adjustwidth}
\begin{adjustwidth}{0.5cm}{0.5cm}
\begin{flushright}
(Second Public Examination: Trinity Term, 1899)
\end{flushright}
\end{adjustwidth}
\vspace{3mm}

Invariant theory, however, did perhaps represent the most sophisticated and modern algebra taught at Oxford in the nineteenth century. It was likely the closest the course came to the abstract algebra that makes up such courses today.

\section*{Number theory}

Questions relating to number theory were a recurring presence on the algebra papers, though the subjects they covered did not show much variation. Modular arithmetic had a large presence in the mathematics course, with related questions appearing consistently throughout the century. These generally asked to show that a number of a specific form was always divisible by a certain number (an example from the first public examination in 1868 is the divisibility of $2n^{3} + 3n^{2} - 5n$ by six for any integer $n$). Wilson's theorem was another especially popular topic. The study of congruences eventually led to the Chinese remainder theorem by the last few decades of the century:

\vspace{3mm}
\begin{adjustwidth}{1cm}{1cm}
\textit{6. Show how to solve the equation $ax - by = c$ in positive integers.}

\textit{Find the general form of all numbers which, when divided by $3, 7, 10$ leave remainders $2, 4, 7$ respectively.}
\end{adjustwidth}
\begin{adjustwidth}{0.5cm}{0.5cm}
\begin{flushright}
(Second Public Examination: Trinity Term, 1890)
\end{flushright}
\end{adjustwidth}
\vspace{3mm}

\noindent
This is another example of something relatively modern making its way into the Oxford course. The Chinese remainder theorem in this non-abstract form had first appeared (without proof) in Europe in 1852, thanks to the missionary Alexander Wylie (Dickson 1920, 57).

Also visible in this question is an example of a linear Diophantine equation. The solvability and number of solutions of such equations came up increasingly in exams over the years. Solvability of equations was of course an important area of research in the nineteenth century, but the only indication of the Oxford course going any further than these simple examples is the following example, which was not followed in later years by anything similar:

\vspace{3mm}
\begin{adjustwidth}{1cm}{1cm}
\textit{3. Shew that neither of the equations $5x^{2} \pm 6y^{2} = z^{2}$ can be solved in rational numbers.}
\end{adjustwidth}
\begin{adjustwidth}{0.5cm}{0.5cm}
\begin{flushright}
(First Public Examination: Trinity Term, 1881)
\end{flushright}
\end{adjustwidth}
\vspace{3mm}

\noindent
Given this, it seems unlikely that undergraduates were exposed to any serious work in this area.

\section*{Algebraic methods in geometry}

The focus of this article is not geometry, but it would be incomplete without mentioning the growing presence algebraic approaches to geometry enjoyed in the nineteenth-century undergraduate course. In the early papers, algebra in the geometry extended as far as expressions for curves to be sketched, with the majority of the rest of the material being inspired by Euclid. This changed with Henry Smith's promotion to the Savilian chair in 1861. Among the subjects he introduced and lectured on was `Modern Geometry', the first such course at any English university (Glaisher 1894, lxxvi). The exam papers indicate that this gave a much more algebraic treatment, and introduced ideas such as projective geometry to undergraduates. There would be questions on these papers that, while geometric in setting, required an almost entirely algebraic approach:

\vspace{3mm}
\begin{adjustwidth}{1cm}{1cm}
\textit{11. Find the equations to the two parabolas which may be drawn through the points of intersection of the conics
\begin{align*}
&7x^{2} + 6xy + 3y^{2} +x +5y - 8 = 0, \\
&3x^{2} + 2xy + y^{2} + x + 2y - 3 = 0.
\end{align*}
}
\end{adjustwidth}
\begin{adjustwidth}{0.5cm}{0.5cm}
\begin{flushright}
(First Public Examination: Michaelmas Term, 1868)
\end{flushright}
\end{adjustwidth}
\vspace{3mm}

The direction of geometry did not change much again for the rest of the century, though there are indications that the idea of certain invariants of algebraic curves had entered the course:

\vspace{3mm}
\begin{adjustwidth}{1cm}{1cm}
\textit{6. Establish Pl\"{u}cker's equations,
\begin{align*}
2\delta + 3\kappa = n(n - 1) - m,\text{ and } \, 6\delta + 8\kappa + \iota = 3n(n - 2)
\end{align*}
for a curve of degree $n$ and class $m$, where $\delta$, $\kappa$, and $\iota$ are respectively the number of nodes, cusps, and points of inflexion.}

\textit{For a cubic curve write down the possible values of $m$, $\delta$, $\kappa$, and $\iota$, and show that a cubic has always one real point of inflexion.}
\end{adjustwidth}
\begin{adjustwidth}{0.5cm}{0.5cm}
\begin{flushright}
(Second Public Examination: Trinity Term, 1899)
\end{flushright}
\end{adjustwidth}
\vspace{3mm}

\section*{Anything more abstract?}

While algebra in the Oxford degree was modernising and there are suggestions of ever more abstract notions there is no evidence that truly abstract algebra was being studied by students at the start of the twentieth century. The course was still very much stuck in the framework of algebraic forms, number theory and geometry without looking at the structure for its own sake. Very occasionally there were questions that hinted at a non-abstract group structure, though without calling it such a thing, as in the following two examples. The first concerns the (cyclic) group of prime roots of unity, and references the idea of a generator. The second, in a more roundabout way, deals with multiplicative inverses in the group of integers modulo $p$.

\vspace{3mm}
\begin{adjustwidth}{1cm}{1cm}
\textit{11. [Show that] All the $n^{th}$ roots of unity can be expressed as the successive powers of any one of the roots except unity, $n$ being a prime number.}

\textit{Solve $x^{7} - 1 = 0$.}
\end{adjustwidth}
\begin{adjustwidth}{0.5cm}{0.5cm}
\begin{flushright}
(First Public Examination: Michaelmas Term, 1875)
\end{flushright}
\end{adjustwidth}
\vspace{3mm}
\begin{adjustwidth}{1cm}{1cm}
\textit{2. [Show that] If $p$ is a prime, the numbers $2, 3, 4, 5, \ldots , p - 2$ can be arranged in $\frac{p-3}{2}$ pairs, such that the product of the members of each pair diminished by unity shall be divisible by $p$.}

\textit{Hence shew that $1 + \lvert \underline{p - 1}$ is divisible by $p$.}
\end{adjustwidth}
\begin{adjustwidth}{0.5cm}{0.5cm}
\begin{flushright}
(Second Public Examination: Michaelmas Term, 1880)\footnote{The second part of the question is Wilson's theorem and irrelevant here; $\lvert \underline{p}$ is obsolete notation for the factorial.}
\end{flushright}
\end{adjustwidth}
\vspace{3mm}

Cayley (1854) made an early attempt at setting out a list of axioms for a group, though his definition was by no means abstract. It was not until the late nineteenth century that something like an abstract group appeared, but the fundamental properties that defined groups had been studied by a handful of mathematicians over the course of the second half of the century (P M Neumann 1999). There would have been ties to these group properties in the theory of linear transformations and invariants, so it is odd that they got no mention in the exam papers at any point. It is perhaps stranger still that the course contained nothing more than hints at the theory of matrices, as this was more similar to invariant theory in that Oxford mathematicians were active in its research. Indeed, matrices were certainly of great interest to Henry Smith, who gave his name to the Smith normal form of a matrix (Smith 1861).

One possible explanation for the slow acceptance of these new ideas into in the Oxford course is a lingering tribalism of Oxford mathematicians. When Cayley, a Cambridge mathematician, was appointed as President for the 1883 meeting in Oxford of the British Association for the Advancement of Science, a group of Oxford scientists and mathematicians organized a petition to protest the choice\footnote{They were unsuccessful in doing anything about the decision, only managing to alienate the association. The BAAS moved the meeting to Southport instead, and didn't return to Oxford until 1894 (Hannabuss 1997, 452).} (Odling 1882). The undersigned included Bartholomew Price, Sedleian Professor of Natural Philosophy, and Charles Pritchard, Savilian Professor of Astronomy. It seems possible then that poor relations between the two universities, and perhaps even bad blood over this specific incident, affected the willingness of Oxford mathematicians to include Cayley's work in their course---invariant theory had of course been taken up by Sylvester by the time it began to appear. This does not offer a fully satisfying explanation for the missing matrices, however, as these were another topic that Sylvester applied himself to. In fact, it is Sylvester himself who is credited with coining `matrix' as a term to describe such an object (Sylvester 1850, 369).

All of this paints a picture of nineteenth-century Oxford as a very insular place, which hampered its own progress where algebra was concerned. Perhaps this is unfair---Hannabuss (1997, 450), for instance, notes that Henry Smith was more highly respected in mainland Europe than in the UK. Certainly, he was awarded the prizes of both the Berlin and Paris Academies of Sciences in 1868 and 1883 respectively. His name was also absent from the aforementioned petition to the BAAS. However, it would appear that despite his own credentials, he was not able to influence the mathematics course to include modern algebra from beyond the Oxford sphere. Nor was Sylvester, who studied for his BA at Cambridge (Parshall 1998, 1) and was a friend of Cayley's. The most likely explanation is perhaps simply that they did not have the clout necessary to push such changes through against opposition from other university figures. In any case, it was left to the twentieth-century Oxford mathematicians to finally introduce abstract algebra to undergraduates.

\newpage

\section*{Acknowledgements}
I would like to thank Dr Christopher Hollings for his valuable guidance and advice in putting this paper together.

\vspace{3mm}
\begin{sloppypar}
\noindent
This is a preprint of an article published by Taylor \& Francis in \textit{BSHM Bulletin: Journal of the British Society for the History of Mathematics} on 10 November 2016, available online: \url{http://www.tandfonline.com/doi/abs/10.1080/17498430.2016.1244630}.
\end{sloppypar}

\section*{Bibliography}

\hangindent=\parindent
\hangafter=1
\noindent
Busbridge, I W, `Oxford mathematics and mathematicians', \textit{University of Oxford mathematical institute website}, Oxford: University of Oxford, \url{https://www.maths.ox.ac.uk/about-us/history/busbridge-lecture}. Accessed on 19 August 2016.

\hangindent=\parindent
\hangafter=1
\noindent
Cayley, A, `On the theory of linear transformations', \textit{Cambridge Mathematical Journal}, 4 (1845), 193--209.

\hangindent=\parindent
\hangafter=1
\noindent
Cayley, A, `On the theory of groups, as depending on the symbolic equation $\theta^{n} = 1$', \textit{Philosophical Magazine}, 7 (1854), 40--47.

\hangindent=\parindent
\hangafter=1
\noindent
Curthoys, M C, `The examination system', in: M G Brock and M C Curthoys (eds), \textit{The history of the University of Oxford}, vol 6, Oxford: Oxford University Press, 1997, 339--74.

\hangindent=\parindent
\hangafter=1
\noindent
Dickson, L E, \textit{History of the theory of numbers}, vol 2, Washington: Carnegie Institution of Washington, 1920.

\hangindent=\parindent
\hangafter=1
\noindent
Glaisher, J W L, Introduction, in: J W L Glaisher (ed), \textit{The collected mathematical papers of Henry John Stephen Smith}, Oxford: Oxford University Press, 1894, lxi--xcv.

\hangindent=\parindent
\hangafter=1
\noindent
Hannabuss, K C, `Mathematics', in: M G Brock and M C Curthoys (eds), \textit{The history of the University of Oxford}, vol 7, Oxford: Oxford University Press, 1997, 443--55.

\hangindent=\parindent
\hangafter=1
\noindent
Hannabuss, K C, `The 19th century', in: J Fauvel, R Flood, and R Wilson (eds), \textit{Oxford figures: eight centuries of the mathematical sciences}, Oxford: Oxford University Press, 2013, 223--38.

\hangindent=\parindent
\hangafter=1
\noindent
Katz, V J, \textit{A history of mathematics: an introduction}, New York: HarperCollins, 1993.

\hangindent=\parindent
\hangafter=1
\noindent
Neumann O, `The Disquisitiones Arithmeticae and the theory of equations', in: C Goldstein, N Schappacher, and J Schwermer (eds), \textit{The shaping of arithmetic after C.F. Gauss's Disquisitiones Arithmeticae}, Berlin: Springer, 2007, 107--28.

\hangindent=\parindent
\hangafter=1
\noindent
Neumann, P M, `What groups were: a study of the development of the axiomatics of group theory', \textit{Bulletin of the Australian Mathematical Society}, 60/2 (1999), 285--301.

\hangindent=\parindent
\hangafter=1
\noindent
Ogilvie, C A, `On the union of classical with mathematical studies', in: D A Talboys (ed), \textit{The Oxford English prize essays}, vol 3, Oxford: Talboys and Browne, 1830, 171--191 (original essay written 1817).

\hangindent=\parindent
\hangafter=1
\noindent
Parshall, K H, \textit{James Joseph Sylvester: life and work in letters}, New York: Oxford University Press, 1998.

\hangindent=\parindent
\hangafter=1
\noindent
Parshall, K H, `Victorian algebra', in: R Flood, A Rice, and R Wilson (eds), \textit{Mathematics in Victorian Britain}, New York: Oxford University Press, 2011, 339--56.

\hangindent=\parindent
\hangafter=1
\noindent
Smith, H J S, `On systems of linear indeterminate equations and congruences', \textit{Philosophical Transactions}, 151 (1861), 293--326.

\hangindent=\parindent
\hangafter=1
\noindent
Sutherland, L S, `The curriculum', in L S Sutherland and L G Mitchell (eds), \textit{The history of the University of Oxford}, vol 5, Oxford: Oxford University Press, 1986, 469--91.

\hangindent=\parindent
\hangafter=1
\noindent
Thompson, S P, \textit{The life of William Thomson}, vol 1, London: Macmillan, 1910.

\hangindent=\parindent
\hangafter=1
\noindent
Weyl, H, `Invariants', \textit{Duke Mathematical Journal}, 5/3 (1939), 489--502. \url{http://projecteuclid.org/euclid.dmj/1077491405}.

\hangindent=\parindent
\hangafter=1
\noindent
Sylvester, J J, `Additions to the articles ``On a new class of theorems'' and ``On Pascal's theorem''', \textit{Philosophical Magazine}, 37 (1850), 363--70.

\subsection*{Unpublished resources}

\hangindent=\parindent
\hangafter=1
\noindent
Odling, W, \textit{Letter to Captain Douglas Galton}, 20 June 1882, Bodleian Library G.A. Oxon 8$^{\circ}$ 1079.

\hangindent=\parindent
\hangafter=1
\noindent
Roberts, B C, \textit{Letters and miscellaneous papers}, 1814, Bodleian Library, 270 c.396.

\hangindent=\parindent
\hangafter=1
\noindent
University of Oxford, \textit{Questions at the mathematical (and physical) examination}, 1830--54, Bodleian Library, 2626 d.21.

\hangindent=\parindent
\hangafter=1
\noindent
University of Oxford, \textit{First public examination in disciplinis mathematicis et physicis}, 1863--72, Bodleian Library, Per. 18753 e.102 (1863/1872).

\hangindent=\parindent
\hangafter=1
\noindent
University of Oxford, \textit{First public examination, mathematics}, 1873--82, Bodleian Library, Per. 18753 e.102.

\hangindent=\parindent
\hangafter=1
\noindent
University of Oxford, \textit{First public examination, mathematics}, 1883--88, Bodleian Library, Per. 18753 e.102.

\hangindent=\parindent
\hangafter=1
\noindent
University of Oxford, \textit{First public examination, mathematics}, 1888--94, Bodleian Library, Per. 18753 e.102.

\hangindent=\parindent
\hangafter=1
\noindent
University of Oxford, \textit{First public examination, mathematics}, 1895--1905, Bodleian Library, Per. 18753 e.102.

\hangindent=\parindent
\hangafter=1
\noindent
University of Oxford, \textit{Second public examination in disciplinis mathematicis et physicis}, 1873--82, Bodleian Library, Per. 18753 e.103.

\hangindent=\parindent
\hangafter=1
\noindent
University of Oxford, \textit{Second public examination in disciplinis mathematicis et physicis}, 1883--88, Bodleian Library, Per. 18753 e.103.

\hangindent=\parindent
\hangafter=1
\noindent
University of Oxford, \textit{Second public examination in disciplinis mathematicis et physicis}, 1888--99, Bodleian Library, Per. 18753 e.103.

\hangindent=\parindent
\hangafter=1
\noindent
University of Oxford, \textit{Second public examination in disciplinis mathematicis et physicis}, 1900--12, Bodleian Library, Per. 18753 e.103.

\end{document}